\documentclass[12pt]{amsart}
\usepackage[letterpaper,margin=1.1in]{geometry}
\usepackage{graphicx}

\newcommand{\RR}{\mathbb{R}}
\newcommand{\BMO}{\mathrm{BMO}}
\newcommand{\VMO}{\mathrm{VMO}}
\newcommand{\Haus}{\mathcal{H}}
\newcommand{\loc}{\mathrm{loc}}
\newcommand{\Tan}{\mathrm{Tan}}
\newcommand{\dist}{\mathrm{dist}}
\newcommand{\spt}{\mathrm{spt}}
\newcommand{\Lip}{\mathrm{Lip}}

\def\res{\hbox{ {\vrule height .22cm}{\leaders\hrule\hskip.2cm} } }

\def\Xint#1{\mathchoice
    {\XXint\displaystyle\textstyle{#1}}%
    {\XXint\textstyle\scriptstyle{#1}}%
    {\XXint\scriptstyle\scriptscriptstyle{#1}}%
    {\XXint\scriptscriptstyle\scriptscriptstyle{#1}}%
    \!\int}
\def\XXint#1#2#3{{\setbox0=\hbox{$#1{#2#3}{\int}$}
    \vcenter{\hbox{$#2#3$}}\kern-.5\wd0}}
\def\dashint{\Xint-}

\newtheorem{theorem}{Theorem}[section]
\newtheorem{prop}[theorem]{Proposition}
\newtheorem{lemma}[theorem]{Lemma}
\newtheorem{cor}[theorem]{Corollary}
\theoremstyle{remark}
\newtheorem{remark}[theorem]{Remark}
\theoremstyle{definition}
\newtheorem{defn}[theorem]{Definition}
\newtheorem{exam}[theorem]{Example}

\numberwithin{equation}{section}

\begin{document}

\title[Harmonic Polynomials and Tangent Measures]{Harmonic Polynomials and Tangent Measures of Harmonic Measure}
\author{Matthew Badger}
\date{October 14, 2009}
\address{Department of Mathematics\\ University of Washington\\ Seattle, WA 98195-4530}
\email{mbadger@math.washington.edu}

\begin{abstract}We show that on an NTA domain if each tangent measure to harmonic measure at a point is a polynomial harmonic measure then the associated polynomials are homogeneous. Geometric information for solutions of a two-phase free boundary problem studied by Kenig and Toro is derived. \end{abstract}
\maketitle

\vspace{-.25cm}

\section{Introduction}

\footnotetext[1]{The author was partially supported by NSF grant DMS-0856687.}

In this paper we use tools from geometric measure theory to catalog fine behavior of harmonic measure on a class of two-sided domains $\Omega\subset\RR^n$, $n\geq 3$. Roughly stated we address the following question. What does a boundary look like if it looks the the same (in terms of harmonic measure) from the interior and from the exterior of a domain? More precisely, if $\Omega$ is 2-sided NTA what conditions does $\partial\Omega$ satisfy when harmonic measure $\omega^+$ on the interior $\Omega^+=\Omega$ and harmonic measure $\omega^-$ on the exterior $\Omega^-=\RR^n\setminus\overline\Omega$ are mutually absolutely continuous? In \cite{KT06}, Kenig and Toro examine this question under the additional hypothesis that the Radon-Nikodym derivative $f=d\omega^-/d\omega^+$ has $\log f\in \VMO(d\omega^+)$. They show that for every point $Q\in \partial\Omega$ and sequence of scales $r_i\downarrow 0$ there is a subsequence (which we relabel) and a harmonic polynomial $h:\RR^n\rightarrow\RR$ such that \begin{equation}\label{firstblowup}\frac{\partial\Omega-Q}{r_i}\rightarrow h^{-1}(0)\quad\text{in Hausdorff distance uniformly on compact sets}.\end{equation} One may hope that only linear polynomials $h$ appear in (\ref{firstblowup}), i.e.\ that the boundary is always flat on small scales; however, there are examples of domains with $\omega^+\ll\omega^-\ll\omega^+$ and $\log f\in C^\infty(\partial\Omega)$ for which non-linear polynomials $h$ appear (see Example 1.4 below). The method in \cite{KT06} relates the geometric blow-ups of the boundary to tangent measures of the harmonic measure. Thus information about the free boundary may be obtained by studying tangent measures of harmonic measure---this is our strategy for the question above. To identify the polynomials appearing in (\ref{firstblowup}), we study properties of ``polynomial harmonic measures" in the topology of weak convergence of Radon measures of $\RR^n$. We prove that only homogeneous harmonic polynomials arise in blow-ups of the boundary.

For any harmonic polynomial $h:\RR^n\rightarrow\RR$, the positive and negative parts $h^\pm$ of $h$ are Green functions with pole at infinity for the unbounded open sets $\{x\in\RR^n:h^\pm(x)>0\}$. The \emph{harmonic measure $\omega_h$ associated to $h$} is the unique harmonic measure with pole at infinity on $\Omega^\pm_h = \{h^\pm>0\}$ with Green function $h^\pm$. That is, for all $\varphi\in C_c^\infty(\RR^n)$, \begin{equation}\label{defnomegah}\int_{\{h=0\}} \varphi\, d\omega_h=\int_{\Omega^+_h} h^+\Delta \varphi = \int_{\Omega^-_h} h^- \Delta \varphi.\end{equation}
Alternatively, by a result of Hardt and Simon \cite{HS}, the zero set $h^{-1}(0)=\partial\Omega_h^\pm$ of a harmonic polynomial is smooth away from a rectifiable subset of Hausdorff dimension at most $n-2$. Hence there exists a unique outward unit normal $\nu^\pm$ on $\partial\Omega_h^\pm$ at almost every point with respect to the surface measure $\sigma=\Haus^{n-1}\res \{h=0\}$ and (\ref{defnomegah}) is equivalent to \begin{equation}d\omega_h=-\frac{\partial h^+}{\partial \nu^+}d\sigma=-\frac{\partial h^-}{\partial\nu^-}d\sigma\end{equation} by the generalized Gauss-Green theorem. In the sequel, we focus on two collections of polynomial harmonic measures that arise as tangent measures of harmonic measure on 2-sided NTA domains examined in \cite{KT06} and \cite{KPT}. (See \S2, \S5 and \S6 below for definitions of tangent measures, NTA and 2-sided NTA domains, respectively.)

Set \begin{equation}\label{defnpd}\mathcal{P}_d=\{\omega_h : h \text{ is a non-zero
harmonic polynomial of degree}\leq d\text{ and }h(0)=0\},\end{equation}\begin{equation}\label{defnfk}\mathcal{F}_k=\{\omega_h: h \text{ is a homogenous harmonic polynomial of degree } k\}.\end{equation} By convention we will use $d$ for the degree of any non-zero polynomial, but reserve $k$ for the degree of a homogeneous polynomial.
If $1\leq k\leq d$, note that $\mathcal{F}_k\subset\mathcal{P}_d$. When $k=1$ the family $\mathcal{F}_1$ is the collection of $(n-1)$-flat measures in $\RR^n$, i.e.\ Hausdorff measures restricted to codimension 1 hyperplanes through the origin.

Our main objective is to exhibit a ``self-improving" property of the tangent measures $\Tan(\omega,Q)$ of harmonic measure $\omega$ at a point $Q$ in the boundary of an NTA domain $\Omega$. Because $\Tan(\omega,Q)$ is independent of the choice of pole for $\omega$ (see Remark \ref{allomgeatan}), we omit the pole from the notation. If $\Omega$ is unbounded, $\omega$ may have a finite pole or pole at infinity.

\begin{theorem}\label{MainThm} Let $\Omega\subset\RR^n$ be a NTA domain with harmonic measure $\omega$. If $Q\in\partial\Omega$ and $\Tan(\omega,Q)\subset\mathcal{P}_d$, then $\Tan(\omega,Q)\subset \mathcal{F}_k$ for some $1\leq k\leq d$.\end{theorem}

The proof of the Theorem \ref{MainThm} illustrates the versatility of a powerful technique from geometric measure theory. Tangent measures are a tool that encode information about the support of a measure, similar to how derivatives describe the local behavior of functions. A remarkable feature is that under general conditions (Theorem \ref{MainGMT}) the cone of tangent measures at a point is connected. This fact lies at the core of Preiss' celebrated paper on rectifiability \cite{P} and recently enabled Kenig, Preiss and Toro \cite{KPT} to compute the Hausdorff dimension of harmonic measure on 2-sided NTA domains with $\omega^+\ll\omega^-\ll\omega^+$. (To appreciate the second result, we invite the reader to compare Theorem 1.2 with the dimension of harmonic measure on Wolff snowflakes \cite{W}, \cite{LVV}.)

\begin{theorem}[\cite{KPT} Theorem 4.3]\label{dimn-1} Let $\Omega\subset\RR^n$ be a 2-sided NTA domain. If harmonic measure $\omega^+$ on the interior $\Omega^+=\Omega$ and harmonic measure $\omega^-$ on the exterior $\Omega^-=\RR^n\setminus\overline{\Omega}$ of $\Omega$ are mutually absolutely continuous, then the Hausdorff dimension of $\omega^\pm$ is $n-1$. Recall this means there exists a subset $\Sigma\subset\partial\Omega$ such that $\dim\Sigma=n-1$ and $\omega^\pm(\partial\Omega\setminus\Sigma)=0$; moreover, if $A\subset\partial\Omega$ and $\dim A<n-1$ then $\omega^\pm(\partial\Omega\setminus A)>0$.
\end{theorem}

In previous instances, connectedness was applied to conclude that the tangent measures of a certain measure (at a.e.\ point) belong to the cone of flat measures $\mathcal{F}_1$. The authors in \cite{KPT} express an opinion that the connectedness of tangent measures ``should be useful in other situations where questions of size and structure of the support of a measure arise." To our knowledge the proof of Theorem 1.1 is the first use of this technique to show that the tangent measures of a measure at a point live in a cone of measures other than $\mathcal{F}_1$.

Stated in the language of tangent measures, Kenig and Toro proved in \cite{KT06} that there exists $d\geq 1$ such that $\Tan(\omega^\pm,Q)\subset\mathcal{P}_d$ for every $Q\in\partial\Omega$. Applying Theorem 1.1 we obtain a refined description of the free boundary. Zooming in along any sequence of scales at a point in the boundary, on a domain satisfying the hypotheses of Theorem 1.3, we see the zero set of a homogeneous harmonic polynomial. The degree of the polynomial is uniquely determined at each point.

\begin{theorem}\label{MainApp} Let $\Omega\subset\RR^n$ be a 2-sided NTA domain with harmonic measure $\omega^+$ on the interior $\Omega^+=\Omega$ and harmonic measure $\omega^-$ on the exterior $\Omega^-=\RR^n\setminus\overline{\Omega}$ of $\Omega$. Assume that $\omega^+$ and $\omega^-$ are mutually absolutely continuous and  $f=d\omega^-/d\omega^+$ satisfies $\log f\in\VMO(d\omega^+)$. Then there exists $d\geq 1$ depending on $n$ and the NTA constants of $\Omega$ and pairwise disjoint sets $\Gamma_1, \dots, \Gamma_d$ such that \begin{equation}\label{gamma1d} \partial\Omega=\Gamma_1\cup\dots\cup\Gamma_d.\end{equation} For each $Q\in \Gamma_k$ and each sequence $r_i\downarrow 0$, there exists a subsequence (which we relabel) and a homogeneous harmonic polynomial $h:\RR^n\rightarrow\RR$ of degree $k$ such that \begin{equation}\frac{\partial\Omega-Q}{r_i}\rightarrow h^{-1}(0)\quad\text{in Hausdorff distance uniformly on compact sets}.\end{equation} Moreover, the domains $\{h^\pm>0\}$ are unbounded 2-sided NTA and  $\omega^\pm(\partial\Omega\setminus\Gamma_1)=0$.\end{theorem}

\begin{exam} In \cite{L} Lewy shows that for $n=3$ there exists a spherical harmonic (homogeneous harmonic polynomial) of degree $k$ whose nodal set divides $S^2$ into two components if and only if $k$ is odd. An explicit example (see Figure 1) is given by
\begin{equation} h(x,y,z)=x^2(y-z)+y^2(z-x)+z^2(x-y)-xyz.\end{equation}
The domain $\Omega=\{h>0\}$ is a 2-sided NTA domain such that for harmonic measures $\omega^+=\omega^-$ with pole at infinity $\log f\equiv 0$ and $0\in\Gamma_3$. Thus, for all $n\geq 3$, it is possible that $\partial\Omega\setminus\Gamma_1$ is non-empty and $\dim\partial\Omega\setminus\Gamma_1\geq n-3$. We do not know if an upper bound on the Hausdorff dimension of $\partial\Omega\setminus\Gamma_1$ holds in general. For instance, is it always true that $\dim\partial\Omega\setminus\Gamma_1<n-1$?

In the plane $(n=2)$ it is known that $\partial\Omega=\Gamma_1$; see Remark 4.3 in \cite{KT06} for details. \hfill$\dashv$\end{exam}
\begin{figure}[h]
\caption{The variety $h^{-1}(0)$ separates the sphere $S^2$ into 2 components.}
\includegraphics[width=2in]{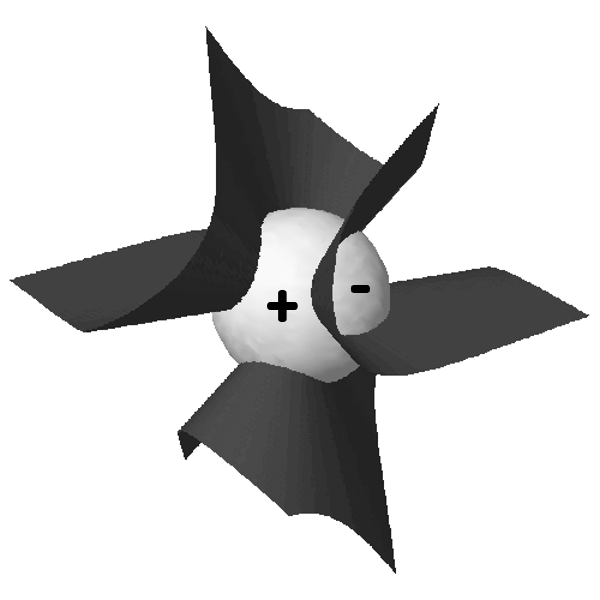} 
\end{figure}

The paper is organized as follows. In \S2 we provide an introduction to tangent measures and related concepts in the general setting of Radon measures on $\RR^n$. The notation established in this section is used pervasively throughout the paper. Our review concludes with an important criterion for connectedness of tangent measures. Here is the rough scheme. Suppose that $\mathcal{M}$ and $\mathcal{F}$ are cones of non-zero Radon measures such that $\mathcal{F}\subset\mathcal{M}$. Furthermore suppose that the set of tangent measures $\Tan(\mu,x)$ of a Radon measure $\mu$ at point $x\in\RR^n$ belongs to $\mathcal{M}$. Under a pair of conditions on $\mathcal{F}$ and $\mathcal{M}$ (see Theorem \ref{MainGMT}) the tangent measures $\Tan(\mu,x)$ are connected relative to $\mathcal{F}$: if one tangent measure of $\mu$ at $x$ belongs to $\mathcal{F}$, then all tangent measures of $\mu$ at $x$ belong to $\mathcal{F}$. While one condition (compactness of $\mathcal{F}$ and $\mathcal{M}$) is routinely checked, verifying the second condition (separation of $\mathcal{F}$ and $\mathcal{M}\setminus\mathcal{F}$) requires work and must be adapted to each situation.

Sections 3 through 5 form the core of the paper. In \S3 we establish inequalities for uniformly bounded spherical harmonics (homogeneous harmonic polynomials restricted to the unit sphere) which depend only on the dimension and degree of the polynomial. In particular, Corollary \ref{bigpiecelemma} is crucial for proving uniform lower estimates for harmonic measures associated to harmonic polynomials of a given degree.

Section 4 studies polynomial harmonic measures in the framework of \S2, focusing on properties which hold independently of assumptions on the underlying domain such as number of components or non-tangential accessibility. The central idea is to consider the rate of doubling at infinity of the measures $\omega_h$, i.e.\ the quantity $\omega_h(B(0,\tau r))/\omega_h(B(0,r))$ as $r\rightarrow\infty$ as a function of $\tau>1$. We show that \begin{equation}\label{doubletheta}\frac{\omega_h(B(0,\tau r))}{\omega_h(B(0,r))}\sim\tau^{n+d-2}\quad\text{as }r\rightarrow\infty, \text{ for every }\tau>1,\end{equation} where $d=\deg h$ and the implied constants for the lower and upper bounds in (\ref{doubletheta}) depend only on $n$ and $d$. Similar bounds for $\omega_h(B(0,\tau r))/\omega_h(B(0,r))$ as $r\rightarrow 0$ are also obtained.

Section 5 is devoted to the proof of Theorem \ref{MainThm}. To start we recall the definition of non-tangentially accessible domains and two useful features of their harmonic measures. The proof of Theorem 1.1 then proceeds in two steps. Suppose that $\Tan(\omega,Q)\subset\mathcal{P}_d$ at some $Q\in\partial\Omega$. Our goal is to show $\Tan(\omega,Q)\subset\mathcal{F}_k$ for some $1\leq k\leq d$. First we apply a blow-up procedure from \cite{KT99} to identify a degree $k=k(Q)$ such that $\Tan(\omega,Q)\cap\mathcal{F}_k\neq\emptyset$. Second we use the doubling property of harmonic measure on NTA domains \cite{JK} and results from section 4 to invoke Theorem \ref{MainGMT} with $\mathcal{F}=\mathcal{F}_k$ and $\mathcal{M}=\Tan(\omega,Q)\cup\mathcal{F}_k$. The connectedness criterion implies that every tangent measure of $\omega$ at $Q$ belongs to $\mathcal{F}_k$.

In \S6, we derive Theorem \ref{MainApp} on homogeneous blow-ups of the boundary of a domain. In addition to Theorem 1.1, we require a blow-up procedure for 2-sided NTA domains from \cite{KT06} and the fact that at almost every point translations of tangent measures are tangent measures. We end by interpreting the decomposition (\ref{gamma1d}) in Theorem \ref{MainApp} from the measure theoretic viewpoint of \S2.

\section{Geometric Measure Theory Ingredients}

Tangent measures and cones of measures were introduced in \cite{P}, where Preiss proved that measures on $\RR^n$ with positive and finite $m$-density almost every are $m$-rectifiable. Here we collect definitions, notation and basic properties of weak convergence of Radon measures, tangent measures and cones of measures which are used throughout the sequel. Much of this material may be found in textbooks of Mattila \cite{M} or Falconer \cite{F}; also see the recent exposition of Preiss' proof by DeLellis \cite{D}. The criterion to check the connectedness of tangent measures (Theorem \ref{MainGMT}) is taken from Kenig-Preiss-Toro \cite{KPT}. Where notations differ across these sources, we adopt the original notation of \cite{P}. (The two novel features of this review are our definition of $F_r$ and the explicit statement of Lemma 2.6.)

Let $B(x,r)$ denote the closed ball with center $x\in\RR^n$ and radius $r>0$. We use the abbreviation $B_r=B(0,r)$ for all $r>0$. Note that $\partial B_1=S^{n-1}$, the unit sphere in $\RR^n$.

A \emph{Radon measure} $\mu$ on $\RR^n$ is a positive Borel regular outer measure on $\RR^n$ that is finite on compact sets. A sequence $(\mu_i)_{i=1}^\infty$ of Radon measures on $\RR^n$ \emph{converges weakly} to a Radon measure $\mu$, written $\mu_i\rightharpoonup \mu$, provided \begin{equation}\label{defnweakconv}\lim_{i\rightarrow\infty} \int fd\mu_i = \int fd\mu\quad\text{for all }f\in C_c(\RR^n).\end{equation} Of course, to test for weak convergence one only needs to check that (\ref{defnweakconv}) holds on a class of functions smaller than $C_c(\RR^n)$; for example, either $C_c^\infty(\RR^n)$ or $\Lip_c(\RR^n)$ suffice. Below we require a quantitative version of weak convergence. To capture the idea that $\mu_i\rightharpoonup\mu$ exactly when $\mu_i$ ``gets close to" $\mu$ on the ball $B_r$ for every (large) $r>0$, we introduce a family of semi-metrics.

Let $\mu$ be a Radon measure on $\RR^n$, and for each $r>0$  define \begin{equation}\label{defnfr}F_r(\mu)=\int_0^r\mu(B_s)ds.\end{equation} Since a Radon measure is locally finite, $F_r(\mu)<\infty$ for all $r>0$. In fact,  \begin{equation}
 \frac{r}2\mu(B_{r/2})\leq F_r(\mu)\leq r\mu(B_r)\quad
 \text{for all }r>0.\end{equation}

If $\mu$ and $\nu$ are Radon measures and $r>0$, we set \begin{equation}F_r(\mu,\nu)=\sup\left\{\left|\int fd\mu-\int fd\nu\right|:f\geq 0, \Lip f\leq 1, \spt f\subset B_r\right\}\end{equation} where $\Lip f$ and $\spt f$ denote the Lipschitz constant and the support of a function $f$, respectively. As an easy exercise one checks $F_r$ is a semi-metric on the set of  Radon measures on $\RR^n$ and a metric on the subset of measures supported in $B_r$. If $r\leq s$ then $F_{r}(\mu,\nu)\leq F_{s}(\mu,\nu)$. Also notice that $F_r(\mu,0)=F_r(\mu)$. Indeed \begin{equation}\label{frequiv}\begin{split} F_r(\mu,0)&=\int\dist(z,\RR^n\setminus B_r)d\mu(z)
=\int_0^\infty \mu\{z:\dist(z,\RR^n\setminus B_r)>s\}ds\\
&=\int_0^r \mu\{z:\dist(z,\RR^n\setminus B_r)>s\}ds=
\int_0^r \mu(B_{r-s})ds=\int_0^r\mu(B_s)ds.
\end{split}\end{equation}

We now state the relationship between weak convergence of Radon measures and $F_r$.

\begin{lemma}[\cite{M} Lemma 14.13]\label{frconv}Suppose that $\mu,\mu_1,\mu_2,\dots$ are Radon measures on $\RR^n$. Then $\mu_i\rightharpoonup\mu$ if and only if $\lim_{i\rightarrow \infty} F_r(\mu_i,\mu)=0$ for all $r>0$.\end{lemma}

\begin{prop}[\cite{P} Proposition 1.12]The Radon measures on $\RR^n$ admit a complete separable metric \begin{equation}\sum_{i=1}^\infty 2^{-i}\min(1,F_i(\mu,\nu))\end{equation} whose topology is equivalent to the topology of weak convergence of Radon measures.\end{prop}

\begin{remark} The family of semi-metrics $F_r$ is related to a distance between probability measures in a compact metric space, which is known by various names in the literature. If $X$ is a compact metric space, the \emph{Kantorovich-Rubinstein} formula \begin{equation}
\sup\left\{\int fd(\mu-\nu):\Lip f\leq 1\right\}\end{equation} defines a complete separable metric on the space of probability measures on $X$ whose topology is equivalent to the weak convergence of probability measures \cite{KR}. For further discussion we refer the reader to the bibliographical notes in Chapter 6 of \cite{V}. \hfill$\dashv$
\end{remark}

Let $x\in\RR^n$ and $r>0$. We write $T_{x,r}$ for the translation by $x$ and dilation by $r$, $T_{x,r}:\RR^n\rightarrow\RR^n$, \begin{equation} T_{x,r}(y)=\frac{y-x}{r}\quad\text{for all }y\in\RR^n.\end{equation} The image measure $T_{x,r}[\mu]$ of a Radon measure $\mu$, which acts on a set $E\subset\RR^n$ by \begin{equation}\label{defntxrmu} T_{x,r}[\mu](E)=\mu(T_{x,r}^{-1}(E))=\mu(x+rE),\end{equation} is also Radon since $T_{x,r}$ is a homeomorphism. In the case $E=B_1$, we interpret (\ref{defntxrmu}) as saying $T_{x,r}[\mu]$ ``blows-up" $B(x,r)$ (for $r$ small) to the unit ball $B_1$ in the sense that $\mu(B(x,r))=T_{x,r}[\mu](B_1)$. Integration against $T_{x,r}[\mu]$ obeys \begin{equation}\int f(z) dT_{x,r}[\mu](z)=\int f\left(\frac{z-x}{r}\right) d\mu(z)\end{equation} whenever at least one of the integrals is defined. Let us pause to record a few simple but highly useful calculations.

\begin{lemma}[Composition Laws] For all $x\in\RR^n$, for all $r,s>0$ and  all measures $\mu,\nu$, \begin{enumerate}
\item $T_{x,rs}=T_{0,s}\circ T_{x,r}$,
\item $T_{x,rs}[\mu]=T_{0,s}[T_{x,r}[\mu]]$,
\item $F_{rs}(\mu)=sF_r(T_{0,s}[\mu])$,
\item $F_{rs}(\mu,\nu)=sF_r(T_{0,s}[\mu],T_{0,s}[\nu])$.
\end{enumerate} \label{compositionlaws}
\end{lemma}

We can now present a definition of tangent measure. The basic idea is to take a sequence of blow-ups $T_{x,r_i}[\mu]$ as $r_i>0$ shrinks to zero and then normalize by some constants $c_i>0$ so that the limit converges.

\begin{defn}Let $\mu$ be a non-zero Radon measure and let $x\in\spt\mu$. We say a non-zero Radon measure $\nu$ is a \emph{tangent measure} of $\mu$ at $x$ and write $\nu\in\Tan(\mu,x)$ if there exists sequences $r_i\downarrow 0$ and $c_i>0$ such that \begin{equation}c_i T_{x,r_i}[\mu]\rightharpoonup \nu.\end{equation}\end{defn}

The set of tangent measures at a point is non-empty under mild assumptions on the measure. For example, if $x\in\spt\mu$ and one of the conditions\begin{itemize}
\item $\overline{D}^s(\mu,x)=\limsup_{r\downarrow 0}\mu(B(x,r))/r^s\in(0,\infty)$ for some $0<s<\infty$
\item $\limsup_{r\downarrow 0}\mu(B(x,2r))/\mu(B(x,r))<\infty$
\end{itemize} hold, then $\Tan(\mu,x)\neq\emptyset$ by the weak compactness of Radon measures.

Taking blow-ups of a measure at a point is closed in the sense that tangent measures to tangent measures are tangent measures. We need two formulations of this principle.

\begin{lemma}\label{tttatzero} Let $\mu$ be a non-zero Radon measure and $x\in\spt\mu$. If $\nu\in\Tan(\mu,x)$, then $\Tan(\nu,0)\subset \Tan(\mu,x)$.\end{lemma}
\begin{proof} Let $\rho\in\Tan(\nu,0)$. Suppose that $r_i,s_i\downarrow 0$ and $c_i,d_i>0$ are sequences such that $c_iT_{x,r_i}[\mu]\rightharpoonup \nu$ and $d_iT_{0,s_i}[\nu]\rightharpoonup \rho$. Since $c_iT_{x,r_i}[\mu]\rightharpoonup\nu$, $\lim_{i\rightarrow\infty} F_1(c_iT_{x,r_i}[\mu],\nu)=0$. Choose a subsequence $(c_{i(j)}, r_{i(j)})$ of $(c_i,r_i)$ such that \begin{equation}F_1(c_{i(j)}T_{x,r_{i(j)}}[\mu],\nu)\leq \frac{1}{j}\left(\frac{s_j}{d_j}\right).\end{equation} After relabeling $(c_{i(j)},r_{i(j)})$, we may assume that \begin{equation}F_1(c_jT_{x,r_j}[\mu],\nu) \leq \frac{1}{j}\left(\frac{s_j}{d_j}\right).\end{equation} Fix $r>0$. Since $F_r$ is a semi-metric, \begin{equation}F_r(c_jd_j T_{x,r_js_j}[\mu],\rho)\leq F_r(c_jd_jT_{x,r_js_j}[\mu],d_jT_{0,s_j}[\nu])
+ F_r(d_jT_{0,s_j}[\nu],\rho).\end{equation} On one hand, $\lim_{j\rightarrow\infty} F_r(d_jT_{0,s_j}[\nu],\rho)=0$, since $d_jT_{0,s_j}[\nu]\rightharpoonup\rho$. On the other hand, for all $j$ sufficiently large such that $s_jr\leq1$, \begin{equation}\begin{split}F_r(c_jd_jT_{x,r_js_j}[\mu],d_jT_{0,s_j}[\nu])
&= d_j F_r(T_{0,s_j}[c_jT_{x,r_j}[\mu]],T_{0,s_j}[\nu])\\
&= \frac{d_j}{s_j}F_{s_jr}(c_jT_{x,r_j}[\mu],\nu)
\leq \frac{d_j}{s_j}F_1(c_jT_{x,r_j}[\mu],\nu)\leq\frac{1}{j}.\end{split}\end{equation}
Hence $\lim_{j\rightarrow\infty} F_r(c_jd_jT_{x,r_js_j}[\mu],\rho)=0$. Since $r>0$ was arbitrary, $c_jd_jT_{x,r_js_j}[\mu]\rightharpoonup \rho$. Therefore, $\Tan(\nu,0)\subset\Tan(\mu,x)$.
\end{proof}

\begin{theorem}[\cite{M} Theorem 14.16]\label{tttae} Let $\mu$ be a non-zero Radon measure. At $\mu$-a.e.\ $x\in\spt\mu$ the following holds: if $\nu\in \Tan(\mu,x)$ and $y\in\spt\nu$, then \begin{enumerate}
\item $T_{y,1}[\nu]\in\Tan(\mu,x)$,
\item $\Tan(\nu,y)\subset\Tan(\mu,x)$.\end{enumerate}\end{theorem}

\begin{proof}[Proof Sketch] The proof of (1) uses the separability of Radon measures in the topology generated by the semi-metrics $F_r$. Statement (2) follows quickly from (1), the composition law $T_{y,r_i}[\nu]=T_{0,r_i}[T_{y,1}[\nu]]$ and Lemma \ref{tttatzero}.\end{proof}

Next we introduce cones of measures or collections of measures which are invariant under scaling.

\begin{defn} A collection $\mathcal{M}$ of non-zero Radon measures on $\RR^n$ is a \emph{cone} provided whenever $\psi\in\mathcal{M}$ and $c>0$ then $c\psi\in\mathcal{M}$. A cone $\mathcal{M}$ is a \emph{d-cone} (or \emph{dilation invariant}) if furthermore $\psi\in\mathcal{M}$ and $r>0$ imply $T_{0,r}[\psi]\in\mathcal{M}$. We also require that $\mathcal{M}\neq\emptyset$.\end{defn}

The technical advantage of working with dilation invariant cones is a simple observation. If $\mathcal{M}$ is a d-cone of Radon measures, then for all $r>0$ there is $\mu\in\mathcal{M}$ such that $F_r(\mu)>0$. Indeed take any $\psi\in\mathcal{M}$. Then $F_{s}(\psi)>0$ for some $s>0$ because $\psi\neq 0$. For any $r>0$, \begin{equation}F_r(T_{0,s/r}[\psi])=\frac{r}{s}F_{r(s/r)}(\psi)=\frac{r}{s} F_{s}(\psi)>0.\end{equation} Since $\mathcal{M}$ is closed under dilations, $\psi_r=T_{0,s/r}[\psi]\in\mathcal{M}$ satisfies $F_r(\psi_r)>0$. In particular, since $F_1(\psi_1)>0$ and $\mathcal{M}$ is closed under scaling the following set is non-empty.

\begin{defn} The \emph{basis} of a d-cone $\mathcal{M}$ is the subset
$\{\psi\in\mathcal{M}:F_1(\psi)=1\}$.\end{defn}

\begin{lemma}[\cite{P} Remark 2.1] Let $\mathcal{M}$ be a d-cone. In the topology of weak convergence of Radon measures, $\mathcal{M}$ is relatively closed (relatively compact) in the collection of all non-zero Radon measures if and only if the basis of $\mathcal{M}$ is closed (compact).\end{lemma}

We are already familiar with the canonical example of a dilation invariant cone.

\begin{lemma}[\cite{P} Remark 2.3] If $\Tan(\mu,x)\neq\emptyset$, then $\Tan(\mu,x)$ is a d-cone with a closed basis.\end{lemma}

Following \cite{P} we define a normalized version of $F_r$ for the distance of a measure to a d-cone of measures as follows. Let $r>0$ and suppose $\sigma$ is a measure such that $F_r(\sigma)>0$. If $\mathcal{M}$ is any d-cone, the ``distance" of $\sigma$ to $\mathcal{M}$ at scale $r$ is given by \begin{equation}d_r(\sigma,\mathcal{M})=
\inf\left\{F_r\left(\frac{\sigma}{F_r(\sigma)},\psi\right):\psi\in\mathcal{M}\text{ and }F_r(\psi)=1\right\}.
\end{equation} If $F_r(\sigma)=0$ we set $d_r(\sigma,\mathcal{M})=1$. Our main use for $d_r$ is to detect, given a pair of nested cones $\mathcal{M}_1\subset\mathcal{M}_2$, if $\mathcal{M}_1$ is separated from $\mathcal{M}_2\setminus\mathcal{M}_1$.

\begin{theorem}[\cite{KPT} Corollary 2.1]\label{MainGMT} Let $\mathcal{F}$ and $\mathcal{M}$ be d-cones, $\mathcal{F}\subset\mathcal{M}$. Assume that \begin{enumerate}
\item Both $\mathcal{F}$ and $\mathcal{M}$ have compact bases,
\item There exists $\epsilon_0>0$ such that whenever $\psi\in\mathcal{M}$ and $d_r(\psi,\mathcal{F})<\epsilon_0$ for all $r\geq r_0$ then $\psi\in\mathcal{F}$.
\end{enumerate} If $\Tan(\mu,x)\subset\mathcal{M}$ and $\Tan(\mu,x)\cap\mathcal{F}\neq\emptyset$, then  $\Tan(\mu,x)\subset\mathcal{F}$.\end{theorem}

We end this review with two conditions that ensure a d-cone has a compact basis. Additional criterion may be found in \cite{P}.

\begin{prop}[\cite{P} Proposition 2.2]\label{closedcompact} Assume $\mathcal{M}$ is a d-cone with a closed basis. Then $\mathcal{M}$ has a compact basis if and only if there exists a finite number $q\geq 1$ such that $\psi(B(0,2r))\leq q\psi(B(0,r))$ for all $\psi\in\mathcal{M}$ and $r>0$.\end{prop}

\begin{cor}[\cite{P} Corollary 2.7]\label{tancompact} Let $\mu$ be a non-zero Radon measure. If $x\in\spt\mu$ and  $\limsup_{r\downarrow 0}\mu(B(x,2r))/\mu(B(x,r))<\infty$ then $\Tan(\mu,x)$ has a compact basis.\end{cor}

\section{Inequalities for Spherical Harmonics}

A well known fact about harmonic functions is that derivatives of a function at a point are controlled by the $L^\infty$-norm of the function in a surrounding ball, in a uniform way depending on the distance of the point to the boundary. Starting from local estimates for the derivatives of harmonic functions on $B_2$ at points of $S^{n-1}$, we derive several inequalities for \emph{spherical harmonics} (homogeneous harmonic polynomials on $\RR^n$ restricted to $S^{n-1}$) of a given degree.

\begin{lemma}\label{dalphabd1} Let $u$ be a real-valued harmonic function on $B_2=\overline{B(0,2)}$. For all $\theta\in S^{n-1}$ and every multi-index $\alpha$, \begin{equation}
|D^\alpha u(\theta)| \leq (2^{n+1}n|\alpha|)^{|\alpha|}\|u\|_{L^\infty(\partial B_2)}.\end{equation}
\end{lemma}

\begin{proof} For example, by Theorem 7 in \S2.2 of \cite{E} with $r=1$, \begin{equation}|D^\alpha u(\theta)|\leq \frac{(2^{n+1}n|\alpha|)^{|\alpha|}}{\omega_n} \|u\|_{L^1(B(\theta,1))}\end{equation} where $\omega_n=\mathcal{L}^n(B(0,1))$ denotes the volume of the unit ball in $\RR^n$. The claim follows since $\|u\|_{L^1(B(\theta,1))}\leq \omega_n \|u\|_{L^\infty(B(\theta,1))}\leq \omega_n\|u\|_{L^\infty(\partial B_2)}$, where the last inequality holds by the maximum principle.
\end{proof}

Uniformly bounded spherical harmonics of degree $k$ have a uniform Lipschitz constant.

\begin{prop}\label{spharmlip} Let $n\geq 2$ and $k\geq 1$. There exists a constant $A_{n,k}>1$ such that for every homogeneous harmonic polynomial $h:\RR^n\rightarrow\RR$ of degree $k$ and every $\theta_1,\theta_2\in S^{n-1}$, \begin{equation}|h(\theta_1)-h(\theta_2)|\leq A_{n,k}\|h\|_{L^\infty(S^{n-1})}|\theta_1-\theta_2|.\end{equation}\end{prop}

\begin{proof} Write $M=\|h\|_{L^\infty(S^{n-1})}$. If $|\theta_1-\theta_2|\geq 1$ then $|h(\theta_1)-h(\theta_2)|\leq 2M\leq 2M|\theta_1-\theta_2|.$

Suppose that $|\theta_1-\theta_2|\leq 1$. By Lemma \ref{dalphabd1}, \begin{equation}\label{dalphabd2}|D^{\alpha} h(\theta_2)| \leq (2^{n+1}n|\alpha|)^{|\alpha|}
\|h\|_{L^\infty(\partial B_2)} = (2^{n+1}n|\alpha|)^{|\alpha|}2^kM\leq(2^{n+2}nk)^kM\end{equation} for every multi-index $\alpha$ with $|\alpha|\leq k$, where $\|h\|_{L^\infty(\partial B_2)}=2^kM$ since $h$ is homogeneous of degree $k$. Expanding $h$ in a Taylor series about $\theta_2$, \begin{equation}\label{taylorh}h(\theta)-h(\theta_2)=\sum_{1\leq|\alpha|\leq k}\frac{D^{\alpha}h(\theta_2)}{\alpha!}(\theta-\theta_2)^\alpha.\end{equation} Evaluating (\ref{taylorh}) at $\theta=\theta_1$ and applying the estimate (\ref{dalphabd2}), \begin{equation} |h(\theta_1)-h(\theta_2)|\leq \sum_{1\leq|\alpha|\leq k} \frac{(2^{n+2}nk)^kM}{\alpha!}|\theta_1-\theta_2|^{|\alpha|}\leq A_{n,k}M|\theta_1-\theta_2|\end{equation} where $A_{n,k}=(2^{n+2}nk)^k\sum_{1\leq|\alpha|\leq k}(\alpha!)^{-1}$.
\end{proof}

The next inequality roughly says that a spherical harmonic takes its ``big values" on a ``big piece" of the unit sphere. Here $\sigma$ denotes surface measure on $S^{n-1}$ with total mass $\sigma(S^{n-1})=\sigma_{n-1}=n\omega_n$.

\begin{cor}\label{bigpiecelemma} Let $n\geq 2$ and $k\geq 1$. There exists a constant $l_{n,k}>0$ such that for every homogeneous harmonic polynomial $h:\RR^n\rightarrow\RR$ of degree $k$, \begin{equation}\sigma\{\theta\in S^{n-1}:|h(\theta)|\geq \tfrac12\|h\|_{L^\infty(S^{n-1})}\}\geq l_{n,k}.\end{equation}\end{cor}

\begin{proof} Choose $\theta_0\in S^{n-1}$ such that $|h(\theta_0)|=\|h\|_{L^\infty(S^{n-1})}=M$. By Proposition \ref{spharmlip}, \begin{equation}
|h(\theta)|\geq |h(\theta_0)|-|h(\theta)-h(\theta_0)|
\geq M(1-A_{n,k}|\theta-\theta_0|).\end{equation}
If $|\theta-\theta_0|\leq 1/2A_{n,k}$, then $|h(\theta)|\geq M/2$. That is, the set $\{\theta\in S^{n-1}:|h(\theta)|\geq M/2\}$ contains the surface ball $\Delta(\theta_0,1/2A_{n,k})$. Thus $l_{n,k}=\sigma(\Delta(\theta_0,1/2A_{n,k}))$ suffices.\end{proof}

Thus the spherical harmonics of degree $k$ satisfy a reverse H\"older inequality.

\begin{cor}\label{reverseholder} Let $n\geq 2$ and $k\geq 1$. There exists a constant $B_{n,k}>1$ such that for every homogeneous harmonic polynomial $h:\RR^n\rightarrow\RR$ of degree $k$, \begin{equation} \|h\|_{L^\infty(S^{n-1})}\leq B_{n,k}\|h\|_{L^1(S^{n-1})}.\end{equation}\end{cor}

\begin{proof} Let $\Gamma=\{\theta\in S^{n-1}:|h(\theta)|\geq \frac12\|h\|_{L^\infty(S^{n-1})}\}$. By Corollary \ref{bigpiecelemma}, \begin{equation}\|h\|_{L^1(S^{n-1})}\geq \frac{1}{2}\|h\|_{L^\infty(S^{n-1})}\sigma(\Gamma) \geq \frac{l_{n,k}}{2}\|h\|_{L^\infty(S^{n-1})}\end{equation} and $B_{n,k}=2/l_{n,k}$ suffices.
\end{proof}

\section{Polynomial Harmonic Measures}

A harmonic polynomial $h:\RR^n\rightarrow\RR$ of degree $d$  decomposes as \begin{equation}\label{hd0}h=h_d+h_{d-1}+\dots+h_0\end{equation} where each non-zero term $h_i$ is a homogenous harmonic polynomial of degree $i$. Indeed if $h=\sum_{|\alpha|\leq d}c_\alpha x^\alpha$ is any polynomial then $h_i=\sum_{|\alpha|=i}c_\alpha x^\alpha$ satisfies (\ref{hd0}). For harmonic $h$, \begin{equation}\label{laphd0}0=\Delta h_d + \Delta h_{d-1} + \dots + \Delta h_{2}.\end{equation} Since $\Delta h_i$ is the sum of monomials of degree $i-2$, the right hand side of (\ref{laphd0}) vanishes only if $\Delta h_i=0$ for all $i\leq d$.

Recall that the collections $\mathcal{P}_d$ and $\mathcal{F}_k$ of polynomial harmonic measures were defined by\begin{itemize}
\item $\mathcal{P}_d=\{\omega_h : h \text{ is a non-zero harmonic polynomial of degree}\leq d\text{ and }h(0)=0\}$,
\item $\mathcal{F}_k=\{\omega_h: h \text{ is a homogenous harmonic polynomial of degree } k\}$.
\end{itemize} Our first observation is that $\mathcal{P}_d$ and $\mathcal{F}_k$ fit into the framework of Section 2.

\begin{lemma} $\mathcal{P}_d$ and $\mathcal{F}_k$ are dilation invariant cones. \label{PdFkAreCones}\end{lemma}
\begin{proof} Suppose that $\omega$ is associated to a harmonic polynomial $h=h_d+\dots+h_1$ and let $c,r>0$. We claim  $cT_{0,r}[\omega]$ is harmonic measure associated to $g(x)=cr^nh(rx)$ where $\Delta g=cr^{n+d}\Delta h_d+\dots+c r^{n+2}\Delta h_2=0$ by the remark following (\ref{laphd0}). For any $\varphi\in C_c^\infty(\RR^n)$, \begin{equation}\begin{split}
\int_{\{g>0\}} g(x)\Delta\varphi(x)dx = \int_{r^{-1}\{h>0\}} cr^nh(rx)\Delta\varphi(x) dx
= c\int_{\{h>0\}} h(y)\Delta\varphi(r^{-1}y)dy\\
= c\int_{\{h=0\}} \varphi(r^{-1}y)d\omega(y)
= c\int_{r^{-1}\{h=0\}} \varphi(x)dT_{0,r}[\omega](x)
= c\int_{\{g=0\}} \varphi(x)dT_{0,r}[\omega](x).\end{split}\end{equation} Since $g$ has the same degree as $h$ and $g$ is homogeneous if $h$ is homogeneous, $\mathcal{P}_d$ and $\mathcal{F}_k$ are dilation invariant cones.\end{proof}

Here is a practical formula to compute $\omega_h$ on balls $B_r$ centered at the origin in terms of the surface measure $\sigma$ on the boundary $\partial B_r$. Throughout this section $\Omega^\pm$ denotes the open sets of positive and negative values of $h$, $\Omega^\pm=\{h^\pm>0\}$.

\begin{lemma}Let $h:\RR^n\rightarrow \RR$ be a harmonic polynomial, $h(0)=0$. For any $r>0$, \begin{equation}\label{l42-1}\omega_h(B_r) = \int_{\partial B_r\cap \Omega^+} \frac{\partial h^+}{\partial r}d\sigma=\int_{\partial B_r\cap \Omega^-}\frac{\partial h^-}{\partial r}d\sigma.\end{equation} If $h$ is homogeneous of degree $k$, then
\begin{equation}\label{l42-2}\omega_h(B_r) =
\frac{k}{2}r^{n+k-2}\|h\|_{L^1(S^{n-1})}.\end{equation}\label{omegaBrAtZero}\end{lemma}

\begin{proof}By a result of Hardt and Simon \cite{HS}, the zero set of a harmonic polynomial is smooth away from a rectifiable subset of dimension at most $n-2$. Hence, for any harmonic polynomial $h:\RR^n\rightarrow\RR$ with $h(0)=0$, the set $B_r\cap\Omega^\pm$ is non-empty and has locally finite perimeter. By the generalized Gauss-Green theorem (c.f.\ Chapter 5 of \cite{E}), \begin{equation} \int_{\partial(B_r\cap\Omega^\pm)}\frac{\partial h^\pm}{\partial\nu^\pm}d\sigma = \int_{B_r\cap \partial\Omega^\pm} \Delta h^\pm = 0\end{equation} where $\nu^\pm$ denotes the unique outer unit normal defined at $\sigma$-a.e.\ $Q\in\partial(B_r\cap \Omega^\pm)$. Thus, writing $\partial(B_r\cap \Omega^\pm)=(\partial B_r\cap\Omega^\pm)\cup(B_r\cap\partial\Omega^\pm)$, \begin{equation}\int_{\partial B_r\cap \Omega^\pm} \frac{\partial h^\pm}{\partial r}d\sigma = -\int_{B_r\cap\partial\Omega^\pm}\frac{\partial h^\pm}{\partial\nu^\pm}d\sigma
=\omega_h(B_r)\end{equation} as desired.

Summing the two formulas in (\ref{l42-1}),
\begin{equation} 2\omega_h(B_r)=\int_{\partial B_r\cap\Omega^+}\frac{\partial h^+}{\partial r}d\sigma + \int_{\partial B_r\cap\Omega^-} \frac{\partial h^-}{\partial r}d\sigma.\end{equation} If $h(r\theta)=r^kh(\theta)$, then $\partial_r h(r\theta)=kr^{k-1}h(\theta)$ and $r\theta\in\Omega^\pm$ if and only if $\theta\in\Omega^\pm$. Hence \begin{equation}\begin{split}
2\omega_h(B_r) &=\int_{\partial B_r\cap\Omega^+}kr^{k-1}h^+(\theta)d\sigma + \int_{\partial B_r\cap\Omega^-}kr^{k-1}h^-(\theta)d\sigma\\&=\int_{\partial B_r}kr^{k-1}|h(\theta)|d\sigma = kr^{n+k-2}\int_{\partial B_1}|h(\theta)|d\sigma\end{split}\end{equation} whenever $h$ is homogeneous of degree $k$.\end{proof}

A consequence of (\ref{l42-2}) is that the measures in $\mathcal{F}_k$ are uniformly doubling at the origin, i.e. for any $\omega\in\mathcal{F}_k$ and $r>0$, \begin{equation}\frac{\omega(B_{2r})}{\omega(B_r)}= 2^{n+k-2}<\infty.\end{equation} We now investigate the doubling properties of measures associated to arbitrary harmonic polynomials. The inequality for spherical harmonics in Corollary \ref{bigpiecelemma} is key.

\begin{lemma}\label{omegaBrNearInfinity} Let $h:\RR^n\rightarrow\RR$ be a harmonic polynomial of degree $d\geq 1$ with $h(0)=0$. There exists $r_1=r_1(n,d,\zeta(h))\geq 1$ such that for all $r>r_1$,\begin{equation} \frac{l_{n,d}}{4}\cdot dr^{n+d-2}\|h_d\|_{L^\infty(S^{n-1})} \leq \omega_h(B_r)\leq \frac{3\sigma_{n-1}}{2} \cdot dr^{n+d-2}\|h_d\|_{L^\infty(S^{n-1})}.\end{equation} Here $\zeta(h)=\max_{1\leq k\leq d-1}\|h_k\|_{L^\infty(S^{n-1})}/\|h_d\|_{L^\infty(S^{n-1})}$
and $r_1=1+12\sigma_{n-1}\zeta(h)/l_{n,d}$.\end{lemma}

\begin{proof} Without loss of generality assume that $M=\|h_d\|_{L^\infty(S^{n-1})}=\|h_d^+\|_{L^\infty(S^{n-1})}$; that is, the maximum of the homogeneous part $h_d$ of $h$ over $S^{n-1}$ is obtained at a positive value. Writing $h$ in polar coordinates,
\begin{align}\label{l43-1}h(r\theta) &= r^d h_d(\theta)+r^{d-1} h_{d-1}(\theta)+\dots + r h_1(\theta),\\
\label{l43-2}\frac{\partial h}{\partial r}(r\theta) &= dr^{d-1}h_d(\theta) +(d-1)r^{d-2}h_{d-1}(\theta)+\dots+h_1(\theta).\end{align} Let $r>1$. Then $\frac1r+\dots+\left(\frac1r\right)^{d-1}\leq \sum_{i=1}^\infty \left(\frac{1}{r}\right)^{i}=\frac{1}{r-1}$ and with $\zeta(h)$ defined as above,
\begin{equation}\label{l43-3}
\left|\frac{r^{d-1}h_{d-1}(\theta)+\dots+rh_1(\theta)}{r^d}\right|
\leq M\zeta(h)\left(\frac{1}{r}+\dots+\frac{1}{r^{d-1}}\right)
\leq \frac{M\zeta(h)}{r-1},\end{equation}
\begin{equation}\label{l43-4} \left|\frac{(d-1)r^{d-2}h_{d-1}(\theta)+\dots+h_1(\theta)}{r^{d-1}}\right|
\leq dM\zeta(h)\left(\frac{1}{r}+\dots+\frac{1}{r^{d-1}}\right)
\leq \frac{d M\zeta(h)}{r-1}.\end{equation} If $r\theta\in \partial B_r\cap\Omega^+$, then $h(r\theta)>0$ and by (\ref{l43-1}) and (\ref{l43-3}),
\begin{equation}\label{l43-5} h_d(\theta) > -\frac{r^{d-1}h_{d-1}(\theta)+\dots+rh_1(\theta)}{r^d} \geq -\frac{M\zeta(h)}{r-1}.\end{equation} Similarly, for all $r>1$ and $\theta\in S^{n-1}$, by (\ref{l43-2}) and (\ref{l43-4}),
\begin{equation}\label{l43-6} dr^{d-1}\left(h_d(\theta)-\frac{M\zeta(h)}{r-1}\right)\leq \frac{\partial h}{\partial r}(r\theta) \leq dr^{d-1}\left(h_d(\theta)+\frac{M\zeta(h)}{r-1}\right).\end{equation}
To estimate $\omega_h(B_r)$ for $r\gg 1$, we will combine (\ref{l42-1}), (\ref{l43-5}) and (\ref{l43-6}) with Corollary \ref{bigpiecelemma}. By the latter, the set $\Gamma=\{\theta\in S^{n-1}:h_d(\theta)\geq M/2\}$ has surface measure $\sigma(\Gamma)\geq l_{n,d}$. Note that $r\Gamma\subset\partial B_r\cap\Omega^+$ provided $r>1+2\zeta(h)$, since $h(r\theta)\geq r^d(h_d(\theta)-M\zeta(h)/(r-1))$, again by (\ref{l43-1}) and (\ref{l43-3}). Put $\Lambda_r=(\partial B_r\cap\Omega^+)\setminus r\Gamma$. Then, by (\ref{l42-1}) and (\ref{l43-6}),\begin{align}
\label{l43-7}\omega_h(B_r)&\geq dr^{d-1}\int_{\partial B_r\cap\Omega^+}\left(h_d(\theta)-\frac{M\zeta(h)}{r-1}\right)d\sigma\\
\label{l43-8}&\geq dr^{d-1}\int_{r\Gamma}\left(\frac{M}{2}-\frac{M\zeta(h)}{r-1}\right)d\sigma+ dr^{d-1}\int_{\Lambda_r}\left(-\frac{M\zeta(h)}{r-1}-\frac{M\zeta(h)}{r-1}\right)d\sigma,
\end{align} where $h_d(\theta)\geq M/2$ on $\Gamma$ by definition and $h_d(\theta)>-M\zeta(h)/(r-1)$ for $r\theta\in\Lambda_r$ by (\ref{l43-5}). Since $\sigma(r\Gamma)\geq l_{n,d}r^{n-1}$ and $\sigma(\Lambda_r)\leq \sigma_{n-1}r^{n-1}$,
\begin{align}
\label{l43-9}\omega_h(B_r)&\geq dr^{d-1}M\left(\frac{1}{2}-\frac{\zeta(h)}{r-1}\right)l_{n,d}r^{n-1}
+dr^{d-1}M\left(-\frac{2\zeta(h)}{r-1}\right)\sigma_{n-1}r^{n-1}\\
\label{l43-10}&\geq dr^{n+d-2}M\left(\frac{l_{n,d}}{2}-\frac{3\sigma_{n-1}\zeta(h)}{r-1}\right).
\end{align} Thus, if $r> 1+12\sigma_{n-1}\zeta(h)/l_{n,d}$, we obtain the lower bound $\omega_h(B_r)\geq (l_{n,d}/4)dr^{n+d-2}M$. A similar (and easier!) estimate using the upper bound in (\ref{l43-6}) shows if $r>1+2\zeta(h)$ then $\omega_h(B_r)\leq (3\sigma_{n-1}/2)dr^{n+d-2}M$. Therefore, it suffices to take $r_1= 1+12\sigma_{n-1}\zeta(h)/l_{n,d}$.
\end{proof}

As an immediate corollary of Lemma \ref{omegaBrNearInfinity} we see that $\omega_h(B_r)$ is doubling as $r\rightarrow\infty$ with doubling constants depending only on $n$ and $d$ in the following sense.

\begin{theorem}\label{doubleatinfinity} There is a constant $C_{n,d}>1$ such that for every $\tau>1$ and every harmonic measure $\omega$ associated to a harmonic polynomial $h:\RR^n\rightarrow\RR$ of degree $d$ with $h(0)=0$, \begin{equation}\frac{\tau^{n+d-2}}{C_{n,d}} \leq \liminf_{r\rightarrow\infty} \frac{\omega(B_{\tau r})}{\omega(B_r)} \leq \limsup_{r\rightarrow\infty} \frac{\omega(B_{\tau r})}{\omega(B_r)} \leq C_{n,d}\tau^{n+d-2}.\end{equation}\end{theorem}

\begin{proof} By Lemma \ref{omegaBrNearInfinity} there exists $r_1\geq 1$ depending on $\omega$ such that for all $r>r_1$, \begin{equation} \frac{l_{n,d}}{6\sigma_{n-1}} \tau^{n+d-2}\leq \frac{\omega(B_{\tau r})}{\omega(B_r)} \leq \frac{6\sigma_{n-1}}{l_{n,d}} \tau^{n+d-2}.\end{equation} Thus, $C_{n,d}=6\sigma_{n-1}/l_{n,d}$ suffices.\end{proof}

While the top degree term of the polynomial $h$ determines the harmonic measure $\omega_h(B_r)$ for large $r$, the non-zero term of lowest degree controls  $\omega_h(B_r)$ on small radii.

\begin{lemma}\label{omegaBrNearZero} Suppose that $h=h_d+h_{d-1}+\dots+h_j$ is a harmonic polynomial with $1\leq j\leq d$ and $h_j\neq 0$. There exists $r_2=r_2(n,j,\zeta_*(h))\leq 1/2$ such that for all $r<r_2$,\begin{equation}\label{l45-0} \frac{l_{n,j}}{4}\cdot jr^{n+j-2}\|h_j\|_{L^\infty(S^{n-1})} \leq \omega_h(B_r)\leq \frac{3\sigma_{n-1}}{2} \cdot jr^{n+j-2}\|h_j\|_{L^\infty(S^{n-1})}.\end{equation} Here $\zeta_*(h)=\max_{j+1\leq k\leq d}\|h_k\|_{L^\infty(S^{n-1})}/\|h_j\|_{L^\infty(S^{n-1})}$
and $r_2=\min(1/2,l_{n,j}/72\sigma_{n-1}\zeta_*(h))$.\end{lemma}

\begin{proof} Without loss of generality assume that $M=\|h_j\|_{L^\infty(S^{n-1})}=\|h_j^+\|_{L^\infty(S^{n-1})}$; that is, the maximum of the homogeneous part $h_j$ of $h$ over $S^{n-1}$ is obtained at a positive value. Writing $h$ in polar coordinates,
\begin{align}\label{l45-1}h(r\theta) &= r^d h_d(\theta)+\dots+r^{j+1} h_{j+1}(\theta)+ r^j h_j(\theta),\\
\label{l45-2}\frac{\partial h}{\partial r}(r\theta) &= dr^{d-1}h_d(\theta) +\dots+(j+1)r^{j}h_{j+1}(\theta)+jr^{j-1}h_j(\theta).\end{align} Let $r\leq 1/2$. Then $r+\dots+r^{d-j}\leq \sum_{i=1}^\infty r^i=\frac{r}{1-r}\leq 2r$ and with $\zeta_*(h)$ defined as above,
\begin{equation}\label{l45-3}
\left|\frac{r^{d}h_{d}(\theta)+\dots+r^{j+1}h_{j+1}(\theta)}{r^j}\right|
\leq M\zeta_*(h)\left(r^{d-j}+\dots+r\right)
\leq 2M\zeta_*(h)r.\end{equation} Also, since $(j+i)/2j\leq i$ for all $i,j\geq 1$ and $\sum_{i=1}^\infty ir^i=\frac{r}{(1-r)^2}\leq 4r$,
\begin{equation}\begin{split}\label{l45-4} \left|\frac{dr^{d-1}h_{d}(\theta)+\dots+(j+1)r^jh_{j+1}(\theta)}{r^{j-1}}\right|
\leq M\zeta_*(h)(dr^{d-j}+\dots+(j+1)r)\\
= 2jM\zeta_*(h)\left(\frac{d}{2j}r^{d-j}+\dots+\frac{j+1}{2j}r\right)
\leq 2jM\zeta_*(h)\sum_{i=1}^\infty ir^i\leq 8jM\zeta_*(h)r.\end{split}\end{equation} If $r\theta\in \partial B_r\cap\Omega^+$, then $h(r\theta)>0$ and by (\ref{l45-1}) and (\ref{l45-3}), \begin{equation}\label{l45-5}h_j(\theta)>-\frac{r^dh_d(\theta)+\dots+r^{j+1}h_{j+1}(\theta)}{r^j}\geq -2M\zeta_*(h)r.\end{equation} Similarly, for all $r\leq 1/2$ and $\theta\in S^{n-1}$, by (\ref{l45-2}) and (\ref{l45-4}),
\begin{equation}\label{l45-6} jr^{j-1}\left(h_j(\theta)-8M\zeta_*(h)r\right)\leq \frac{\partial h}{\partial r}(r\theta) \leq jr^{j-1}\left(h_j(\theta)+8M\zeta_*(h)r\right).\end{equation} By Corollary 3.3, the set $\Gamma=\{\theta\in S^{n-1}:h_j(\theta)\geq M/2\}$ has surface measure $\sigma(\Gamma)\geq l_{n,j}$. Note $r\Gamma\subset\partial B_r\cap \Omega^+$ if $r<1/4\zeta_*(h)$, since $h(r\theta)\geq r^j(h_j(\theta)-2M\zeta_*(h)r)$, again by (\ref{l45-1}) and (\ref{l45-3}). Put $\Lambda_r=(\partial B_r\cap\Omega^+)\setminus r\Gamma$. Then, by (\ref{l42-1}) and (\ref{l45-6}), \begin{align}
\label{l45-7}\omega_h(B_r)&\geq jr^{j-1}\int_{\partial B_r\cap\Omega^+}\left(h_j(\theta)-8M\zeta_*(h)r\right)d\sigma\\
\label{l45-8}&\geq jr^{j-1}M\int_{r\Gamma}\left(\frac{1}{2}-8\zeta_*(h)r\right)d\sigma+ jr^{j-1}M\int_{\Lambda_r}\left(-2\zeta_*(h)r-8\zeta_*(h)r\right)d\sigma,
\end{align} where $h_j(\theta)\geq M/2$ on $\Gamma$ by definition and $h_j(\theta)>-2M\zeta_*(h)r$ for $r\theta\in\Lambda_r$ by (\ref{l45-5}). Since $\sigma(r\Gamma)\geq l_{n,j}r^{n-1}$ and $\sigma(\Lambda_r)\leq \sigma_{n-1}r^{n-1}$, if $r<1/16\zeta_*(h)$ we obtain
\begin{align}
\label{l45-9}\omega_h(B_r)&\geq jr^{j-1}M\left(\frac{1}{2}-8\zeta_*(h)r\right)l_{n,j}r^{n-1}
+jr^{j-1}M\left(-10\zeta_*(h)r\right)\sigma_{n-1}r^{n-1}\\
\label{l45-10}&\geq jr^{n+j-2}M\left(\frac{l_{n,j}}{2}-18\sigma_{n-1}\zeta_*(h)r\right).
\end{align} Thus, if $r<\min(1/2, l_{n,j}/72\sigma_{n-1}\zeta_*(h))$, we get the lower bound \begin{equation}\omega_h(B_r)\geq (l_{n,j}/4)jr^{n+j-2}M.\end{equation} The estimate $\omega_h(B_r)\leq (3\sigma_{n-1}/2)jr^{n+j-2}M$ for all $r< \min(1/2,1/16\zeta_*(h))$ follows easily from (\ref{l42-1}) and the upper bound in $(\ref{l45-6})$. Therefore, the estimates (\ref{l45-0}) for $\omega_h(B_r)$ hold for all $r<r_2$ with $r_2=\min(1/2,l_{n,j}/72\sigma_{n-1}\zeta_*(h))$.
\end{proof}

\begin{theorem}\label{doubleatzero} There is a constant $c_{n,j}>1$ such that for every $\tau>1$ and every harmonic measure $\omega$ associated to a polynomial $h=h_d+h_{d-1}+\dots+h_j$ with $1\leq j\leq d$ and $h_j\neq 0$,
\begin{equation}\label{doubleatzeroeq}\frac{\tau^{n+j-2}}{c_{n,j}} \leq \liminf_{r\rightarrow 0} \frac{\omega(B_{\tau r})}{\omega(B_r)} \leq \limsup_{r\rightarrow 0} \frac{\omega(B_{\tau r})}{\omega(B_r)} \leq c_{n,j}\tau^{n+j-2}.\end{equation}
\end{theorem}
\begin{proof} By Lemma \ref{omegaBrNearZero} there exists $r_2\leq 1/2$ depending on $\omega$ such that whenever $\tau r<r_2$, \begin{equation} \frac{l_{n,j}}{6\sigma_{n-1}} \tau^{n+j-2}\leq \frac{\omega(B_{\tau r})}{\omega(B_r)} \leq \frac{6\sigma_{n-1}}{l_{n,j}} \tau^{n+j-2}.\end{equation} Thus, $c_{n,j}=6\sigma_{n-1}/l_{n,j}$ suffices.
\end{proof}

The next lemma generalizes Lemma 4.1 in \cite{KPT}; notice that the assumption $\{h>0\}$ and $\{h<0\}$ are NTA domains has been removed.

\begin{lemma}\label{epsilon0}Suppose $h:\RR^n\rightarrow\RR$ is a harmonic polynomial of degree $d\geq 1$ with $h(0)=0$, and let $\omega$ be harmonic measure associated to $h$. There exists $\epsilon_0>0$ depending only on $n$, $d$ and $k$ such that if $d_r(\omega,\mathcal{F}_k)<\epsilon_0$ for all $r\geq r_0$ then $d=k$. \end{lemma}

\begin{proof} Let $\tau>1$ and choose $r\geq r_0$ such that $d_{\tau r}(\omega,\mathcal{F}_k)<\epsilon_0$. Then there exists  $\psi\in\mathcal{F}_k$ such that $F_{\tau r}(\psi)=1$ and
\begin{equation}\label{l46-1}F_r\left(\frac{\omega}{F_{\tau r}(\omega)},\psi\right)\leq
F_{\tau r}\left(\frac{\omega}{F_{\tau r}(\omega)},\psi\right)<\epsilon_0.\end{equation} Hence, by the triangle inequality,
\begin{equation}\label{l46-2}F_r(\psi)-\epsilon_0<\frac{F_r(\omega)}{F_{\tau r}(\omega)}<F_r(\psi)+\epsilon_0.\end{equation} Since $\psi$ is associated to a homogeneous polynomial of degree $k$, say $p$, by Lemma \ref{omegaBrAtZero},
\begin{equation}\label{l46-3}F_r(\psi)=\int_0^r \psi(B_s)ds=\frac{k \|p\|_{L^1(S^{n-1})}}{2}\int_0^r s^{n+k-2}ds = \frac{k \|p\|_{L^1(S^{n-1})}}{2(n+k-1)}r^{n+k-1}\end{equation} for all $r>0$.
In particular, $1=F_{\tau r}(\psi)=\tau^{n+k-1} F_r(\psi)$. That is,
\begin{equation}\label{l46-4} F_r(\psi)=\tau^{-n-k+1}.\end{equation} Moreover, since $(r/2)\omega(B_{r/2})\leq F_r(\omega) \leq r \omega(B_r)$ for all $r$, by Theorem  \ref{doubleatinfinity},
\begin{equation}\label{l46-5}\frac{1}{C_{n,d}}\left(\frac{1}{2\tau}\right)^{n+d-1}\leq \frac{1}{2\tau}\frac{\omega(B_{r/2})}{\omega(B_{\tau r})}\leq \frac{F_r(\omega)}{F_{\tau r}(\omega)}\leq \frac{2}{\tau}\frac{\omega(B_r)}{\omega(B_{\tau r/2})} \leq C_{n,d}\left(\frac{2}{\tau}\right)^{n+d-1}\end{equation} for all $r>r_1(h)$. Setting $\widetilde C=C_{n,d}2^{n+d-1}>1$, \begin{equation}\label{l46-6}\widetilde{C}^{-1}\tau^{-n-d+1}\leq \frac{F_r(\omega)}{F_{\tau r}(\omega)}\leq \widetilde C \tau ^{-n-d+1}\end{equation} for all $r> r_1$. Combining (\ref{l46-2}), (\ref{l46-4}) and (\ref{l46-6}) yields
\begin{equation}\label{l46-7}\tau^{-n-k+1}-\epsilon_0 < \widetilde C \tau^{-n-d+1}\quad\text{and}\quad \widetilde C^{-1}\tau^{-n-d+1}< \tau^{-n-k+1}+\epsilon_0.\end{equation} Equivalently, \begin{equation}\label{l46-8}\tau^{d-k}(1-\tau^{n+k-1}\epsilon_0)<\widetilde C\quad\text{and}\quad  \tau^{k-d}(1+\tau^{n+k-1}\epsilon_0)^{-1}<\widetilde C.\end{equation}
Because $\widetilde C$ is independent of $\tau$, we can set $\tau =2\widetilde C$. Thus, for $(2\widetilde C)^{n+k-1}\epsilon_0= 1/2$,
\begin{equation}\label{l46-9} \frac{1}{2}(2\widetilde C)^{d-k}<\widetilde C\quad\text{and}\quad \frac{2}{3}(2\widetilde C)^{k-d}<\widetilde C.\end{equation} On a moment's reflection one sees (\ref{l46-9}) is impossible if $d\neq k$. (For example, if $d-k\geq 1$, then $\widetilde C=\frac{1}{2}(2\widetilde C)
\leq \frac{1}{2}(2\widetilde C)^{d-k}<\widetilde C$. If $k-d\geq 1$, then $\frac43\widetilde C=\frac23(2\widetilde C)\leq \frac23(2\widetilde C)^{k-d}<\widetilde C$.) Therefore, if $d_r(\omega,\mathcal{F}_k)<\epsilon_0=\frac12(2\widetilde C)^{-n-k+1}$ for all $r\geq r_0$ then $h$ has degree $k$.\end{proof}

For emphasis let us remark again that $\epsilon_0$ in Lemma \ref{epsilon0} only depends on the dimension, the degree $d$ of the polynomial $h$ and the degree $k$ of the ``homogeneous cone" $\mathcal{F}_k$. Taking the minimum of  finitely many $\epsilon_0$ from Lemma \ref{epsilon0} we obtain:

\begin{cor}\label{epsilon1}There is $\epsilon_1=\epsilon_1(n,d)>0$ with the property if $\omega\in\mathcal{P}_d$ and $d_r(\omega,\mathcal{F}_k)<\epsilon_1$ for all $r\geq r_0$ with $1\leq k\leq d$ then the degree of the polynomial associated to $\omega$ is $k$.\end{cor}

\begin{cor}\label{epsilon2}There is $\epsilon_2=\epsilon_2(n,d)>0$ with the property if $\omega\in\mathcal{P}_d$ and $d_r(\omega,\mathcal{F}_1)<\epsilon_2$ for all $r\geq r_0$ then $\omega\in \mathcal{F}_1$.\end{cor}

In order to invoke Theorem \ref{MainGMT} the cones studied must satisfy a compactness condition. Recall that the basis of a dilation invariant cone $\mathcal{M}$ is $\{\psi\in\mathcal{M}:F_1(\psi)=1\}$.

\begin{lemma}For each $k\geq 1$, $\mathcal{F}_k$ has a compact basis. \label{FkCompact}\end{lemma}

\begin{proof} First we claim there exists a constant $C=C(n,k)<\infty$ such that the coefficients of any polynomial associated to a harmonic measure in the basis of $\mathcal{F}_k$ are bounded by $C$. Let $\omega\in\mathcal{F}_k$ satisfying $F_1(\omega)=1$ be associated to the homogeneous harmonic polynomial $h$ of degree $k$. By (\ref{l42-2}) and the definition of $F_1$, \begin{equation}\label{l49-1}
F_1(\omega)=\int_0^1\omega(B_s)ds=\frac{k}{2(n+k-1)}\|h\|_{L^1(S^{n-1})}.
\end{equation} Since $F_1(\omega)=1$, $\|h\|_{L^1(S^{n-1})}=2(n+k-1)/k$. Hence, by Corollary \ref{reverseholder},
\begin{equation}\label{l49-2}\|h\|_{L^\infty(S^{n-1})}\leq B_{n,k} \|h\|_{L^1(S^{n-1})}= \frac{2B_{n,k}(n+k-1)}{k}.\end{equation} If $h(X)=\sum_{|\alpha|=k}c_\alpha X^\alpha$ then $|c_\alpha|=|D^\alpha h(0)|/\alpha!\leq
|D^{\alpha}h(0)|$ by Taylor's formula. Then the mean value property for $D^\alpha h$ and estimate (\ref{dalphabd2}) yield
\begin{equation}\label{l49-3} |c_{\alpha}|\leq \dashint_{S^{n-1}}|D^\alpha h(\theta)|d\sigma(\theta)
\leq \sup_{\theta\in S^{n-1}}|D^\alpha h(\theta)|\leq (2^{n+2}nk)^k\|h\|_{L^\infty(S^{n-1})}.\end{equation} Combining (\ref{l49-2}) and (\ref{l49-3}) shows that $|c_\alpha|\leq C(n,k)$ for every coefficient of $h$.

Now let $\omega^i\in\mathcal{F}_k$ be any sequence of measures such that $F_1(\omega^i)=1$, and let $h^i$ be the polynomial associated to $\omega^i$. By the argument above, the coefficients of $h^i$ are uniformly bounded. Hence from $h^i$ we can extract a subsequence $h^{i_j}\rightarrow h^\infty$ uniformly on compact subsets of $\RR^n$, where $h^\infty$ is either identically zero or a homogeneous harmonic polynomial of degree $k$. (We will exclude the first possibility shortly). If $\varphi\in C_c^\infty(\RR^n)$, then \begin{equation}\lim_{j\rightarrow\infty} \int\varphi d\omega^{i_j}=\lim_{j\rightarrow\infty}\int (h^{i_j})^+\Delta\varphi=\int (h^\infty)^+\Delta\varphi=\int\varphi d\omega_{h^\infty}.\end{equation} Thus $\omega^{i_j}\rightharpoonup \omega^\infty=\omega_{h^\infty}$ and since $F_1(\omega^\infty)=\lim_{j\rightarrow\infty} F_1(\omega^{i_j})=1$, $h^\infty\not\equiv 0$. We have shown that for every sequence  $\omega^i\in\mathcal{F}_k$ with $F_1(\omega^i)=1$ there is a subsequence $\omega^{i_j}\rightharpoonup\omega^\infty\in\mathcal{F}_k$. Therefore, $\mathcal{F}_k$ has a compact basis.
\end{proof}

We do not know if the cone $\mathcal{P}_d$ has a closed or compact basis for $d\geq 2$. To implement the method of Lemma \ref{FkCompact} and show that $\mathcal{P}_d$ has a compact basis, one must find a way to control $\|h\|_{L^\infty(S^{n-1})}$ from the data $F_1(\omega_h)=1$. On the other hand, to prove that $\mathcal{P}_d$ does not have a compact basis, by Proposition $\ref{closedcompact}$ it suffice to produce a sequence of measures $\omega_i\in\mathcal{P}_d$ and radii $r_i>0$ such that $\sup_i\omega_i(B_{2r_i})/\omega_i(B_{r_i})=\infty$. Since polynomial harmonic measures are doubling near infinity (Theorem \ref{doubleatinfinity}) and doubling near zero (Theorem \ref{doubleatzero}), candidate radii must be selected from an intermediate range depending on $\zeta(h)$ and $\zeta_*(h)$. The main challenge lies in estimating $\omega_h(B_r)$ on these middle scales. Since $\zeta(h)\zeta_*(h)\leq 1$ for every quadratic polynomial $h$, the final answer may depend on whether $d=2$ or $d\geq 3$.

\section{Polynomial Tangent Measures are Homogeneous}

We now recast our focus to polynomial harmonic measures which appear as tangent measures of harmonic measure on NTA domains and take up the proof of Theorem 1.1. Jerison and Kenig introduced non-tangentially accessible domains in $\RR^n$ as a natural class of domains on which Fatou type convergence theorems hold for harmonic functions \cite{JK}. Here the doubling of harmonic measure on NTA domains is combined with properties from Section 4 and a blow-up procedure from \cite{KT99} in order to invoke Theorem \ref{MainGMT}.

We start by recalling the definitions of NTA domains.

\begin{defn}An open set $\Omega\subset\RR^n$ satisfies the \emph{corkscrew condition} with constants $M>1$ and $R>0$ provided for every $Q\in\partial\Omega$ and $0<r<R$ there exists a \emph{non-tangential point} $A=A(Q,r)\in\Omega$ such that $M^{-1}r<|A-Q|<r$ and $\dist(A,\partial\Omega)> M^{-1}r$.\end{defn}

An $M$--\emph{non-tangential ball} $B(X,r)$ in a domain $\Omega$, is an open ball contained in $\Omega$ whose distance to $\partial\Omega$ is comparable to its radius in the sense that \begin{equation}
M^{-1}r<\dist(B(X,r),\partial\Omega)<Mr.\end{equation} For $X_1,X_2\in\Omega$ a \emph{Harnack chain} form $X_1$ to $X_2$ is a sequence of $M$--non-tangential balls such that the first ball contains $X_1$, the last  contains $X_2$, and consecutive balls intersect.

\begin{defn}A domain $\Omega\subset\RR^n$ satisfies the \emph{Harnack chain condition} with constants $M>1$ and $R>0$ if for every $Q\in\partial\Omega$ and $0<r<R$ when $X_1,X_2\in \Omega\cap B(Q,\frac{r}{4})$ satisfy \begin{equation} \min_{j=1,2}\dist(X_j,\partial\Omega)>\varepsilon\quad\text{and}\quad|X_1-X_2|< 2^k\varepsilon\end{equation} then there is a Harnack chain from $X_1$ to $X_2$ of length $Mk$ such that the diameter of each ball is bounded below by $M^{-1}\min_{j=1,2}\dist(X_j,\partial\Omega)$.\end{defn}

\begin{defn}A domain $\Omega\subset\RR^n$ is \emph{non-tangentially accessible} or \emph{NTA} if there exist $M>1$ and $R>0$ such that (i) $\Omega$ satisfies the corkscrew and Harnack chain conditions, (ii) $\RR^n\setminus\overline{\Omega}$ satisfies the corkscrew condition. If $\partial\Omega$ is unbounded then we require $R=\infty$.\end{defn}

A bounded simply connected domain $\Omega\subset\RR^2$ is NTA if and only if $\Omega$ is a quasidisk (the image of the unit disk under a quasiconformal map of the plane). In higher dimensions, while every quasiball (the image of the unit ball under a quasiconformal map of $\RR^n$, $n\geq 3$) is still a bounded NTA domain, there exist bounded NTA domains homeomorphic to a ball in $\RR^n$ which are not quasispheres. The reader may consult \cite{JK} for more information. Also see \cite{KT97} where it is shown that every $\delta$-Reifenberg flat domain in $\RR^n$ with $\delta<\delta_n$ is non-tangentially accessible.

Harmonic measure on NTA domains is locally doubling. While Jerison and Kenig only considered bounded domains, their proof of this result extends to the unbounded case.

\begin{lemma}[\cite{JK} Lemmas 4.8, 4.11]\label{omegaFiniteDoubling} Let $\Omega\subset\RR^n$ be a NTA domain. There exists a constant $C<\infty$ depending on the NTA constants of $\Omega$ such that if $Q\in\partial\Omega$, $0<2r<R$ and $X\in\Omega\setminus B(Q,2Mr)$ then $\omega^X(B(Q,2s))\leq C\omega^X(B(Q,s))$ for all $0<s<r$.\end{lemma}

On an unbounded NTA domain there is a related doubling measure called harmonic measure with pole at infinity, which is obtained as the weak limit of harmonic measures $\omega^{X_i}$ (properly rescaled) as $X_i\rightarrow\infty$.

\begin{lemma}[\cite{KT99} Lemma 3.7, Corollary 3.2]\label{omegaInfiniteDoubling} Let $\Omega\subset\RR^n$ be an unbounded NTA domain. There exists a doubling Radon measure $\omega^\infty$ supported on $\partial\Omega$ satisfying \begin{equation} \int_{\partial\Omega} \varphi\, d\omega^\infty = \int_\Omega u\Delta\varphi\quad\text{for all }\varphi\in C_c^\infty(\RR^n)\end{equation} where \begin{equation} \left\{\begin{array}{rl} \Delta u=0 &\text{in }\Omega\\ u>0 &\text{in }\Omega\\ u=0 &\text{on }\partial\Omega.\end{array}\right.\end{equation} The measure $\omega^\infty$ and Green function $u$ are unique up to multiplication by a positive scalar. We call $\omega^\infty$  a harmonic measure of $\Omega$ with pole at infinity.\end{lemma}

When a result about harmonic measure of a domain $\Omega$ is independent of the choice of pole, we denote the measure by $\omega$ without any superscript. This means that when $\Omega$ is unbounded we allow $\omega$ to have a finite pole or pole at infinity.

\begin{lemma}\label{ntatancompact} If $\Omega\subset\RR^n$ is NTA and $Q\in\partial\Omega$, then $\Tan(\omega,Q)$ has a compact basis.\end{lemma}

\begin{proof} At any point in the support, the tangent measures of an asymptotically doubling measure has a compact basis by Corollary \ref{tancompact}. This is true on an NTA domain by Lemma \ref{omegaFiniteDoubling} when $\omega$ has a finite pole and by Lemma \ref{omegaInfiniteDoubling} when $\omega$ has pole at infinity.\end{proof}

On an NTA domain there is a correspondence between the tangent measures of harmonic measure and geometric blow-ups of the domain and boundary \cite{KT99}. Let $\Omega\subset\RR^n$ be a NTA domain, let $Q\in\partial\Omega$ and let $r_i\downarrow 0$. For each $i$, zoom in on the domain, the boundary and the harmonic measure at $Q$ and scale $r_i$: \begin{equation}\label{1blowupdefn}
\Omega_i = \frac{\Omega-Q}{r_i},\quad \partial\Omega_i=\frac{\partial\Omega-Q}{r_i},\quad \omega_i=\frac{T_{Q,r_i}[\omega]}{\omega(B(Q,r_i))}.\end{equation}

\begin{theorem}[\cite{KT99} Lemma 3.8]\label{ktblowup1} Let $\Omega\subset\RR^n$ be a NTA domain with harmonic measure $\omega$, let $Q\in\partial\Omega$ and let $r_i\downarrow 0$. Define $\Omega_i$, $\partial\Omega_i$ and $\omega_i$ by (\ref{1blowupdefn}). There exists a subsequence of $r_i$ (which we relabel) and an unbounded NTA domain $\Omega_\infty\subset\RR^n$ such that \begin{equation} \Omega_i\rightarrow \Omega_\infty\quad\text{in Hausdorff distance sense uniformly on compact sets},\end{equation}
\begin{equation} \partial\Omega_i\rightarrow \partial\Omega_\infty\quad\text{in Hausdorff distance sense uniformly on compact sets}.\end{equation} Moreover,
\begin{equation} \omega_i\rightharpoonup \omega_\infty\end{equation} where $\omega_\infty$ is harmonic measure for $\Omega_\infty$ with pole at infinity.\end{theorem}

\begin{remark}\label{allomgeatan} The measure $\omega_\infty$ in Theorem \ref{ktblowup1} obtained as a weak limit of the blow-ups $\omega(B(Q,r_i))^{-1}T_{Q,r_i}[\omega]$ is a tangent measure of $\omega$ at $Q$. In fact, up to scaling by a constant, every tangent measure of $\omega$ at $Q$ has this form since $\omega$ is doubling; c.f.\ \cite{M} Remark 14.4. Hence, since the blow-ups $\Omega_i$ of the domain $\Omega$ do not depend on the pole of harmonic measure, the cone of tangent measures $\Tan(\omega,Q)$ is also independent of the pole of $\omega$. \hfill$\dashv$
\end{remark}

The next lemma identifies the degree $k$ of the cone $\mathcal{F}_k$ appearing in Theorem 1.1.

\begin{lemma}\label{choosemin}Let $\Omega\subset\RR^n$ be a NTA domain, let $Q\in\partial\Omega$, and assume $\Tan(\omega,Q)\subset\mathcal{P}_{d}$. If $k$ is the minimum degree such that $\mathcal{P}_k\cap\Tan(\omega,Q)\neq\emptyset$, then $\mathcal{P}_k\cap\Tan(\omega,Q)\subset\mathcal{F}_k$.\end{lemma}

\begin{proof} If $k=1$, then $\mathcal{P}_1=\mathcal{F}_1$. If $k\geq 2$, suppose for contradiction that there exists $\nu\in\Tan(\omega,Q)$ associated to a nonhomogeneous harmonic polynomial $h$ of degree $k$, say $h=h_k+h_{k-1}+\dots+h_j$ with $j<k$ and $h_j\neq 0$. By Theorem \ref{ktblowup1} (applied to $\Omega$ and $\omega$), either  $\{x\in\RR^n:h(x)>0\}$ or $\{x\in\RR^n:h(x)<0\}$ is an unbounded NTA domain where $\nu$ is a harmonic measure with pole at infinity for that domain. Without loss of generality, assume $U=\{x\in\RR^n:h(x)>0\}$ is an unbounded NTA domain and $\nu$ is harmonic measure on $U$ with pole at infinity. Choose $r_i\downarrow 0$. By Theorem \ref{ktblowup1} (now applied to $U$ and $\nu$), there is a subsequence $r_i$ and an unbounded NTA domain $U_\infty$ such that \begin{equation}U_i = \frac{U}{r_i}\rightarrow U_\infty\quad\text{and}\quad \partial U_i = \frac{\partial U}{r_i}\rightarrow\partial U_\infty\end{equation} in the sense of Hausdorff distance uniformly on compact sets and \begin{equation}\nu_i=\frac{T_{0,r_i}[\nu]}{\nu(B_{r_i})} \rightharpoonup\nu_\infty.\end{equation} Moreover, $\nu_\infty$ is harmonic measure with pole at infinity for $U_\infty$. Observe that $\partial U_i$ is the set of all $y\in\RR^n$ such that $h(r_iy)=0$, i.e. \begin{equation}r_i^{k}h_k(y)+r_i^{k-1}h_{k-1}(y)+ \dots+r_i^jh_j(y)=0.\end{equation} Dividing by $r_i^j$ and letting $i\rightarrow\infty$, we see $\partial U_\infty$ is the set of all $y\in\RR^n$ such that $h_j(y)=0$ and $\nu_\infty\in\mathcal{P}_j$. By Lemma \ref{tttatzero}, $\nu_{\infty}\in\Tan(\omega,Q)\cap\mathcal{P}_j$ is a blow up of $\omega$ corresponding to a harmonic polynomial of degree $j<k$. This contradicts the minimality of $k$. Therefore, every blow up of $\omega$ at $Q$ of minimum degree is homogeneous, i.e.  $\mathcal{P}_k\cap\Tan(\omega,Q)\subset\mathcal{F}_k$.
\end{proof}

We now have all the pieces to prove Theorem 1.1. Recall: \emph{Let $\Omega\subset\RR^n$ be a NTA domain with harmonic measure $\omega$. If $Q\in\partial\Omega$ and $\Tan(\omega,Q)\subset\mathcal{P}_{d}$, then $\Tan(\omega,Q)\subset\mathcal{F}_k$ for some $1\leq k\leq d$.}

\begin{proof}[Proof of Theorem 1.1] Let $k=\min\{j:\mathcal{P}_j\cap\Tan(\omega,Q)\neq \emptyset\}\leq d$ and set \begin{equation}\mathcal{F}=\mathcal{F}_k,\quad \mathcal{M}=\Tan(\omega,Q)\cup\mathcal{F}_k.\end{equation} Then $\mathcal{F}\subset\mathcal{M}$ and both d-cones have a compact basis by Lemma \ref{FkCompact} and Lemma \ref{ntatancompact}. Since $\mathcal{M}\subset\mathcal{P}_d$, Corollary \ref{epsilon1} and Lemma \ref{choosemin} together imply that there exists an $\epsilon_1>0$ such that for all $\mu\in\mathcal{M}$ if $d_r(\mu,\mathcal{F}_k)<\epsilon_1$ for all $r\geq r_0$ then $\mu\in\mathcal{F}_k$. By Theorem \ref{MainGMT} (the connectedness of tangent measures), since $\Tan(\omega,Q)\subset\mathcal{M}$ and $\Tan(\omega,Q)\cap\mathcal{F}_k\neq \emptyset$, we conclude $\Tan(\omega,Q)\subset\mathcal{F}_k$.
\end{proof}

\section{Blow-ups on 2-sided NTA Domains}

\begin{defn}A domain $\Omega\subset\RR^n$ is \emph{two-sided non-tangentially accessible} or \emph{2-sided NTA} if $\Omega^+=\Omega$ and $\Omega^-=\RR^n\setminus\overline{\Omega}$ are NTA; i.e., there are $M>1$ and $R>0$  such that $\Omega^\pm$ satisfy the corkscrew and Harnack chain conditions, and if $\partial\Omega$ is unbounded we require  $R=\infty$.\end{defn}

Throughout this section we use the convention that if $\Omega\subset\RR^n$ is a 2-sided domain, then $\omega^+$ is harmonic measure on the interior $\Omega^+=\Omega$ and $\omega^-$ is harmonic measure on the exterior $\Omega^-=\RR^n\setminus\overline{\Omega}$ of $\Omega$. If $\Omega^+$ or $\Omega^-$ is unbounded, then we allow $\omega^+$ or $\omega^-$ to have a finite pole or pole at infinity, respectively.

There is a two-sided version of the blow-up procedure for NTA domains \cite{KT06}. Let $\Omega\subset\RR^n$ be a 2-sided NTA domain, let $Q\in\partial\Omega$ and let $r_i\downarrow 0$. Let $u^\pm$ be the Green function for $\Omega^\pm$ with the same pole as the harmonic measure $\omega^\pm$. We zoom in on the interior and exterior domains, boundary, harmonic measures and Green functions at $Q$ along scales $r_i$:
\begin{equation}\label{2blowupdefn}
\Omega^\pm_i=\frac{\Omega^\pm-Q}{r_i},\quad
\partial\Omega_i = \frac{\partial\Omega-Q}{r_i},\quad
\omega^\pm_i = \frac{T_{Q,r_i}[\omega^\pm]}{\omega^\pm(B(Q,r_i))},\quad
u^\pm_i = \frac{u^\pm \circ T_{Q,r_i}^{-1}}{\omega^\pm(B(Q,r_i))}r_i^{n-2}.\end{equation}

\begin{theorem}[\cite{KT06} Theorem 4.2]\label{ktblowup2} Let $\Omega\subset\RR^n$ be a 2-sided NTA domain, $Q\in\partial\Omega$ and  $r_i\downarrow 0$.  Define the sets $\Omega_i^\pm$ and $\partial\Omega_i$, measures $\omega_i^\pm$ and functions $u_i^\pm$ by (\ref{2blowupdefn}). There is a subsequence of $r_i$ (which we relabel) and an unbounded 2-sided NTA domain $\Omega_\infty$ such that \begin{equation}\Omega_i^\pm\rightarrow\Omega_\infty^\pm\quad\text{in Hausdorff distance uniformly on compact sets},\end{equation}
\begin{equation}\partial\Omega_i\rightarrow\partial\Omega_\infty\quad\text{in Hausdorff distance uniformly on compact sets}.\end{equation} Moreover,
\begin{equation}\omega_i^\pm\rightharpoonup\omega^\pm_\infty,\end{equation} \begin{equation}u_i^\pm\rightarrow u_\infty^\pm\quad\text{uniformly on compact sets}\end{equation} where $\omega^\pm_\infty$ is harmonic measure with pole at infinity for $\Omega^\pm$ and Green function $u_\infty^\pm$.\end{theorem}

Two more definitions are necessary.

\begin{defn}Let $\Omega\subset\RR^n$ be a NTA domain with harmonic measure $\omega$. We say that $f\in L^2_\loc(d\omega)$ has \emph{bounded mean oscillation} with respect to $\omega$, i.e. $f\in\BMO(d\omega)$ if \begin{equation}\sup_{r>0}\sup_{Q\in\partial\Omega}\left(\dashint_{B(Q,r)} |f-f_{Q,r}|^2d\omega\right)^{1/2}<\infty\end{equation} where $f_{Q,r}=\dashint_{B(Q,r)}fd\omega$.\end{defn}

\begin{defn}\label{defnVMO}Let $\Omega\subset\RR^n$ be a NTA domain with harmonic measure $\omega$. Let $\VMO(d\omega)$ denote the closure of the set of bounded uniformly continuous functions defined on $\partial\Omega$ in $\BMO(d\omega)$. If $f\in\VMO(d\omega)$ we say $f$ has \emph{vanishing mean oscillation}.\end{defn}

Polynomial harmonic measures appear as tangent measures on domains with mutually absolutely continuous interior and exterior harmonic measures.

\begin{theorem}\label{kt44}Let $\Omega\subset\RR^n$ be a 2-sided NTA domain with interior harmonic measure $\omega^+$ and exterior harmonic measure $\omega^-$. Assume that $\omega^+\ll\omega^-\ll\omega^+$ and $f=d\omega^-/d\omega^+$ satisfies  $\log f\in\VMO(d\omega^+)$. There exists $d\geq 1$ depending on $n$ and the NTA constants of $\Omega$ such that $\Tan(\omega^+,Q)=\Tan(\omega^-,Q)\subset\mathcal{P}_d$ for all $Q\in\partial\Omega$.\end{theorem}

\begin{proof} Under the same hypothesis, Theorem 4.4 in \cite{KT06} concludes that, in the notation of Theorem 6.2 above, $\omega^+_\infty=\omega^-_\infty$ and $u=u_\infty^+-u_\infty^-$ is a harmonic polynomial. The proof that $u$ is a polynomial shows there exists $d\geq 1$ determined by $n$ and the NTA constants of $\Omega$ such that $u$ has degree at most $d$. The correspondence between tangent measures and the blow-ups of Green functions in Theorem \ref{ktblowup2} implies $\Tan(\omega^+,Q)=\Tan(\omega^-,Q)\subset\mathcal{P}_d$.\end{proof}

The self-improving property of tangent measures in Theorem 1.1 yields:

\begin{cor}\label{MainCor} Let $\Omega\subset\RR^n$ be a $2$-sided NTA domain with interior harmonic measure $\omega^+$ and exterior harmonic measure $\omega^-$. Assume that $\omega^+\ll\omega^-\ll\omega^+$ and $f=d\omega^-/d\omega^+$ satisfies $\log f\in \VMO(d\omega^+)$. There exists $d\geq 1$ depending on $n$ and the NTA constants of $\Omega$ and pairwise disjoint sets $\Gamma_1, \dots, \Gamma_d$ such that \begin{equation}\partial\Omega = \Gamma_1\cup\dots\cup\Gamma_{d},\end{equation} where $\Tan(\omega^+,Q)=\Tan(\omega^-,Q)\subset\mathcal{F}_k$ for all $1\leq k\leq d$ and $Q\in\Gamma_k$.\end{cor}

The decomposition of the boundary in Corollary \ref{MainCor} has an extra interpretation from the geometric measure theory viewpoint. Unfortunately the proof of Theorem \ref{tttae} does not provide a certificate to check at which points in the support of a measure the translations of tangent measures are tangent measures. But the corollary identifies the points in the support of harmonic measure where this behavior occurs. To state the result, we first write down a precise definition of the desired property.

\begin{defn}\label{omegaTTT} Let $\mathcal{M}$ be a cone of non-zero Radon measures on $\RR^n$. We say that $\mathcal{M}$ is \emph{translation invariant} if $T_{x,1}[\mu]\in\mathcal{M}$ for all $\mu\in\mathcal{M}$ and all $x\in\spt\mu$. \end{defn}

\begin{prop} \label{tttdecomp} Let $\Omega$ be as in Corollary \ref{MainCor}. Then the cone $\Tan(\omega^\pm,Q)$ is translation invariant if and only if $Q\in\Gamma_1$.\end{prop}

\begin{proof} If $\mu$ is a flat measure, then $T_{x,1}[\mu]=\mu$ for every $x\in\spt\mu$. Hence $\Tan(\omega^\pm,Q)\subset\mathcal{F}_1$ is translation invariant for every $Q\in\Gamma_1$.

Conversely, assume $\Tan(\omega^\pm,Q)\subset\mathcal{F}_k$ is translation invariant and let $\nu\in\Tan(\omega^\pm,Q)$. Then $\spt\nu=h^{-1}(0)$ for some harmonic polynomial $h$. By \cite{HS} the zero set of a harmonic polynomial is smooth away from a rectifiable subset of dimension at most $n-2$. Hence, $\spt\nu$ is smooth at some $x\in\spt\nu$. Because $T_{x,1}[\nu]\in\Tan(\omega^\pm,Q)$ and $\spt T_{x,1}[\nu]=\spt\nu-x$, we conclude there exists $\sigma\in\Tan(\omega^\pm,Q)\subset\mathcal{F}_k$ such that $\spt\sigma$ is smooth at $0$. But the zero set of a non-zero homogeneous polynomial of degree $k$ (the support of $\sigma=T_{x,1}[\nu]$) is smooth at $0$ only if $k=1$. Therefore, $\Tan(\omega^\pm,Q)\subset\mathcal{F}_1$ and $Q\in\Gamma_1$. \end{proof}

\begin{cor}\label{g2dnull} Let $\Omega$ be as in Corollary \ref{MainCor}. If $d\geq 2$, then $\omega^\pm(\Gamma_2\cup\dots\cup\Gamma_d)=0$.\end{cor}

\begin{proof} By Theorem \ref{tttae}, the cone $\Tan(\omega^\pm,Q)$ of tangent measures at $Q$ is translation invariant for $\omega^\pm$-a.e.\ $Q\in\partial\Omega$. Since this property fails at all $Q\in\Gamma_2\cup\dots\cup\Gamma_d$, the set must have zero harmonic measure.\end{proof}

We can now record:

\begin{proof}[Proof of Theorem 1.3] Let $\Omega\subset\RR^n$ be a 2-sided NTA domain such that $\omega^+\ll\omega^-\ll\omega^+$ and $\log d\omega^-/d\omega^+\in \VMO(d\omega^+)$. By Corollary \ref{MainCor} we can write $\partial\Omega=\Gamma_1\cup\dots\cup\Gamma_d$ where $\Tan(\omega^\pm,Q)\subset\mathcal{F}_k$ for all $Q\in\Gamma_k$ (and $d$ only depends on $n$ and the NTA constants of $\Omega$). By Corollary \ref{g2dnull}, $\omega^\pm(\partial\Omega\setminus\Gamma_1)=\omega^\pm(\Gamma_2\cup\dots\cup\Gamma_d)=0$.

Suppose that $Q\in\Gamma_k$ and we are given $r_i\downarrow 0$. By Theorem \ref{ktblowup2}, there is a subsequence of $r_i$ (which we relabel) such that $\omega^+(B(Q,r_i))^{-1}T_{Q,r_i}[\omega^+]\rightharpoonup \omega^+_\infty\in\Tan(\omega^+,Q)$ and \begin{equation}\frac{\partial\Omega-Q}{r_i}\rightarrow \spt\omega^+_\infty\quad\text{in Hausdorff distance uniformly on compact sets}.
\end{equation} Since $\omega^+_\infty\in\mathcal{F}_k$, there exists a homogeneous harmonic polynomial $h:\RR^n\rightarrow\RR$ of degree $k$ such that $\spt\omega^+_\infty=h^{-1}(0)$.\end{proof}

\begin{remark} One can also apply Theorem 1.1 to tangent measures on two-sided domains without any assumptions on the Radon-Nikodym derivative $d\omega^-/d\omega^+$. Let $\Omega\subset\RR^n$ be an arbitrary 2-sided NTA domain. First we recall the definition of the set $\Gamma\subset\partial\Omega$ from \cite{KPT}. By the differentiation theory of Radon measures, \begin{equation}h(Q)=\lim_{r\downarrow 0}\frac{\omega^-(B(Q,r))}{\omega^+(B(Q,r))}\in[0,\infty]\end{equation} exists at $\omega^\pm$-a.e.\ $Q\in\partial\Omega$. Let \begin{equation}\Lambda=\{Q\in\partial\Omega:h(Q)\text{ exists, } 0<h(Q)<\infty\}.\end{equation} It is easily seen that $\omega^+\ll\omega^-\ll\omega^+$ on $\Lambda$ and $\omega^+\perp\omega^-$ on $\partial\Omega\setminus\Lambda$. (Note that \cite{KPT} uses the notation `$\Lambda_1$' for $\Lambda$. They also define sets $\Lambda_2$, $\Lambda_3$ and $\Lambda_4$ which we do not need here.) To define $\Gamma$ we restrict our attention to density points of $\Lambda$ and $h$: \begin{equation}\Gamma=\{Q\in\Lambda: Q\text{ is a density point of } \Lambda \text{ and a Lebesgue point of }h\text{ w.r.t. }\omega^+\}.\end{equation} Note $\Gamma$ agrees with $\Lambda$ up to a set of $\omega^\pm$ measure zero and any subset $A\subset\partial\Omega$ such that $\omega^+\res A\ll\omega^-\res A\ll \omega^+\res A$ can be written as $A=B\cup N$ where $\omega^\pm(N)=0$ and $B\subset\Gamma$. Thus, up to a set of $\omega^\pm$ measure zero, $\Gamma$ is the maximal ``mutually absolutely continuous" piece of $\partial\Omega$. By Theorems 3.3 and 3.4 in \cite{KPT} (analogously to Theorems \ref{ktblowup2} and \ref{kt44}), there exists $d\geq 1$ such that $\Tan(\omega^+,Q)=\Tan(\omega^-,Q)\subset\mathcal{P}_{d}$ if $Q\in\Gamma$. Hence, by Theorem 1.1, \begin{equation}\Gamma=\Gamma_1\cup\dots\cup\Gamma_d,\end{equation} where for each $Q\in\Gamma_k$, $\Tan(\omega^+,Q)=\Tan(\omega^-,Q)\subset\mathcal{F}_k$.

In particular, if $\Omega\subset\RR^n$ is a 2-sided NTA domain and $\omega^+\ll\omega^-\ll\omega^+$, then \begin{equation}\partial\Omega=\Gamma\cup N=\Gamma_1\cup\dots\cup\Gamma_d\cup N\end{equation} where $\omega^\pm(N)=0$ and $\Tan(\omega^+,Q)=\Tan(\omega^-,Q)\subset\mathcal{F}_k$ for each $Q\in\Gamma_k$. \hfill$\dashv$
\end{remark}

\end{document}